\documentclass[12pt]{amsart}

\usepackage{euscript}
\usepackage{amsmath}
\usepackage{amsthm}
\usepackage{epsfig}
\usepackage{amssymb}

\usepackage{epic}

\DeclareMathOperator{\coker}{{\rm coker}}

\numberwithin{equation}{section}

\theoremstyle{plain}
\newtheorem{theorem}[equation]{Theorem}
\newtheorem{lemma}[equation]{Lemma}
\newtheorem{proposition}[equation]{Proposition}

\newtheorem{conj/thm}[equation]{Conjecture/Theorem}

\theoremstyle{definition}
\newtheorem{example}[equation]{Example}
\newtheorem{remark}[equation]{Remark}

\newtheorem{definition}[equation]{Definition}

\title{On the log canonical threshold and numerical data of a resolution in dimension 2}

\author{Willem Veys}
\address{KU Leuven, Dept. Wiskunde, Celestijnenlaan 200B, 3001 Leuven, Belgium}
\email{wim.veys@kuleuven.be
}

\thanks{The  author is partially supported by
 KU Leuven grant C14/17/083 and wants to thank Kien Nguyen for discussing the two--dimensional case.}


\date{}

\begin{document}

\begin{abstract}
We show various properties of numerical data of an embedded resolution of singularities for plane curves, which are inspired by a conjecture of Igusa on exponential sums.
\end{abstract}

\maketitle


\pagestyle{myheadings} \markboth{{\normalsize  W.
Veys}}{ {\normalsize Numerical data}}

\renewcommand{\div}{{\rm div}}

\newcommand{\et}{\mathcal{T}}
\newcommand{\bS}{{\mathbb S}}
\newcommand{\bma}{\mbox{\boldmath$a$}}
\newcommand{\bmb}{\mbox{\boldmath$b$}}
\newcommand{\bmc}{\mbox{\boldmath$c$}}
\newcommand{\bme}{\mbox{\boldmath$e$}}
\newcommand{\bmi}{\mbox{\boldmath$i$}}
\newcommand{\bmj}{\mbox{\boldmath$j$}}
\newcommand{\bmv}{\mbox{\boldmath$v$}}
\newcommand{\bmk}{\mbox{\boldmath$k$}}
\newcommand{\bmm}{\mbox{\boldmath$m$}}
\newcommand{\bms}{\mbox{\boldmath$s$}}
\newcommand{\bSW}{\mbox{\boldmath$SW$}}
\newcommand{\bmf}{\mbox{\boldmath$f$}}
\newcommand{\bmg}{\mbox{\boldmath$g$}}
\newcommand{\gq}{{\mathfrak q}}
\newcommand{\xo}{o}
\newcommand{\veeK}{{\tiny\vee}}
\newcommand{\gh}{g}

\newcommand{\lp}{{l}}
\newcommand{\ev}{\varepsilon}
\newcommand{\tx}{\tilde{X}}
\newcommand{\tz}{\tilde{Z}}
\newcommand{\calL}{{\mathcal L}}
\newcommand{\calm}{{\mathcal M}}
\newcommand{\calx}{{\mathcal X}}
\newcommand{\calo}{{\mathcal O}}
\newcommand{\calt}{{\mathcal T}}
\newcommand{\cali}{{\mathcal I}}
\newcommand{\calj}{{\mathcal J}}
\newcommand{\calC}{{\mathcal C}}
\newcommand{\calS}{{\mathcal S}}
\newcommand{\calQ}{{\mathcal Q}}
\newcommand{\calF}{{\mathcal F}}

\newcommand{\cs}{\langle \chi_0\rangle}

\newcommand{\cc}{\bar{C}}
\newcommand{\vp}{\varphi}

\def\mmod{\mbox{mod}}
\let\d\partial
\def\EE{\mathcal E}
\newcommand{\cC}{\EuScript{C}}
\def\C{\mathbb C}
\def\Q{\mathbb Q}
\def\R{\mathbb R}
\def\bS{\mathbb S}
\def\bH{\mathbb H}
\def\bB{\mathbb B}\def\bC{\mathbb C}\def\bA{\mathbb A}
\def\Z{\mathbb Z}
\def\N{\mathbb N}
\def\bn{\mathbb N}
\def\bp{\mathbb P}\def\bt{\mathbb T}
\def\eop{$\hfill\square$}
\def\bif{(\, , \,)}
\def\coker{\mbox{coker}}
\def\im{{\rm Im}}

\newcommand{\Gammma}{{G}}
\newcommand{\no}{\noindent}
\newcommand{\bfc}{{\mathbb C}}
\newcommand{\bfq}{{\mathbb Q}}
\newcommand{\calE}{{\mathcal E}}
\newcommand{\calW}{{\mathcal W}}
\newcommand{\calV}{{\mathcal V}}
\newcommand{\calP}{{\mathcal P}}
\newcommand{\calI}{{\mathcal I}}\newcommand{\calJ}{{\mathcal J}}
\newcommand{\calA}{{\mathcal A}}\newcommand{\CalA}{{\calA_F\cup\calA_W}}
\newcommand{\calAA}{{\mathcal A}'}
\newcommand{\calB}{{\mathcal B}}
\newcommand{\calR}{{\mathcal R}}
\newcommand{\calG}{{\mathcal G}}\newcommand{\calN}{{\mathcal N}}
\newcommand{\bc}{{\mathbb C}}
\newcommand{\bez}{B_{\epsilon_0}}
\newcommand{\br}{{\mathbb R}}
\newcommand{\bq}{{\mathbb Q}}
\newcommand{\sez}{S_{\epsilon_0}}
\newcommand{\ep}{\epsilon}
\newcommand{\vs}{\vspace{3mm}}
\newcommand{\si}{\sigma}
\newcommand{\Gammas}{S}
\newcommand{\Gx}{G_\pi(X)}
\newcommand{\Gxf}{G_\pi(X,f)}
\newcommand{\GxF}{G_\pi(X,F)}
\newcommand{\Gax}{\Gamma_\pi(X)}
\newcommand{\Gaxf}{\Gamma_\pi(X,f)}
\newcommand{\GaxF}{\Gamma_\pi(X,F)}
\newcommand{\GaxFW}{\Gamma_\pi(X,F,W)}

\newcommand{\q}{w}

\newcommand{\labelpar}{\label}



\section{Introduction} Singularity invariants of a hypersurface are often described in terms of a chosen embedded resolution. In particular, the so--called
numerical data  $(N_i, \nu_i)$ of a resolution are crucial in various invariants, e.g. poles of zeta functions of Igusa type \cite{I2}\cite{DL1}\cite{DL2}, jumping coefficients of multiplier ideals \cite{ELSV}, roots of Bernstein--Sato polynomials \cite{K}, monodromy eigenvalues \cite{AC},  etc.
In particular $\min_i \frac {\nu_i}{N_i}$ does not depend on the chosen resolution, and is nowadays called the log canonical threshold, see e.g. \cite{M}.

Let $f$ be a polynomial in $n$ variables. In a previous version of the manuscript \cite{CMN} an equivalence was shown between a statement on numerical data of an embedded resolution of $f$ and a famous old conjecture on exponential sums of Igusa \cite{I1} (as well as with a local version of that conjecture by Denef--Sperber \cite{DS} and a more general version by Cluckers--Veys \cite{CV}).

In the present paper, we present a proof of that statement in dimension $n=2$. After finishing this work, we learned that Cluckers--Musta\c t\u a--Nguyen proved the statement in arbitrary dimension, using techniques from the Minimal Model Program. We think however that various aspects of our more elementary proof of the two--dimensional case are of independent interest.

In \S 2 we fix notation and state the Conjecture/Theorem on numerical data of Cluckers--Musta\c t\u a--Nguyen. An important ingredient in our proof for $n=2$ is a new property/formula for the numbers $\nu_i$ in terms of the dual resolution graph for plane curves, with appropriate decorations along edges, which we establish in \S 3. Then, in \S 4, we show (a somewhat stronger version of) the statement in \cite{CMN}.

\section{Preliminaries}

Let $f \in \C[x]=\C[x_1,\dots,x_n]$ be a nonconstant polynomial. Fix an embedded resolution $\pi$ of $f$, that is, $\pi:Y\to \C^n$ is a proper birational morphism satisfying (i) $Y$ is a (complex) nonsingular algebraic variety, (ii) $\pi$ is an isomorphism outside $\pi^{-1}\{f=0\}$, (iii) $\pi^{-1}\{f=0\}$ is a simple normal crossings divisor.

We denote by $E_i, i\in T,$ the (nonsingular) irreducible components of $\pi^{-1}\{f=0\}$. Let $N_i$ and $\nu_i-1$ denote the multiplicity of $E_i$ in the divisor of $\pi^*f=f\circ \pi$ and $\pi^*(dx_1\wedge\dots\wedge dx_n)$, respectively. In other words, $\div(f \circ \pi) = \sum_{i\in T} N_iE_i$ and the canonical divisor $K_Y=K_{\pi}=\sum_{i\in T} (\nu_i-1)E_i$.
The $(N_i,\nu_i)_{i\in T}$ are called the {\em numerical data} of $\pi$.

In order to formulate the statement of \cite{CMN}, we need the notion of power condition. The normal crossings condition says that, for any point $P\in Y$, there is an affine neighbourhood $V$ of $P$, such that
\begin{equation}\label{nc}
 f\circ\pi = u\prod_{i\in I} y_i^{N_i} ,
\end{equation}
for some $I\subset T$, in the coordinate ring $\mathcal{O}_V$ of $V$. Here $i\in I$ if and only if $P\in E_i$, $u$ is a unit in $\mathcal{O}_V$, the component $E_i$ is given by $y_i=0$ and the $(y_i)_{i\in I}$ form a regular sequence in the local ring of $Y$ at $P$.

Here we only state the local version of the power condition; this is the relevant one for the present paper. Also, we only need the local version of the log canonical threshold.  Assuming that $f(0)=0$, the {\em log canonical threshold of $f$ at $0\, (\in \C^n)$} is
$$
c_0 = \min_{i\in T, 0\in \pi(E_i)} \nu_i/N_i .
$$

\begin{definition}\label{power}
 Let $f$ and $\pi$ be as above. Let $d\in \Z_{\geq 2}$.
We say that $(f,\pi)$ {\em satisfies the $d$--power condition}
if there exists a nonempty open $W$ in some irreducible component of some  $\cap_{i\in I} E_i$ with $\pi(W)=\{0\}$  and some $g\in  \mathcal{O}_W$ such that
$$
d\mid N_i \text{ for all } i\in I
$$
and
$$
u |_W = g^d ,
$$
where $u$ is as in (\ref{nc}) on an open $V$ satisfying $W= (\cap_{i\in I} E_i) \cap V$.
\end{definition}

\begin{conj/thm}[\cite{CMN}]\label{CNconj}
Let $f$ and $\pi$ be as above. If $(f,\pi)$ satisfies the $d$--power condition, then
\begin{equation}\label{CNinequality}
c_0 \leq \frac 1d + \sum_{i\in I} N_i (\frac{\nu_i}{N_i} - c_0 ).
\end{equation}
\end{conj/thm}

We prefer to rewrite this inequality in the form
\begin{equation}\label{formulation}
 c_0 \leq \frac{\sum_{i\in I} \nu_i + 1/d}{\sum_{i\in I} N_i + 1} .
\end{equation}

\smallskip
\begin{remark}\label{trivial}
We note that Conjecture \ref{CNconj} is trivial when $c_0 \leq 1/d$. This happens in particular when $d \mid N_\ell$ for some component $E_\ell$ of the strict transform of $f$.
\end{remark}

When $|I|=n$ in Definition \ref{power}, the power condition is automatically satisfied, since then $W=\cap_{i\in I} E_i$ is a point and $u |_W$ is a constant. Otherwise, the following property is a useful corollary of the power condition. It is shown in \cite{CMN}, but follows also from an easy local computation, using unique factorization in a regular local ring.

\begin{proposition}\label{d-divisible}
Let $f$ and $\pi$ be as above. If $(f,\pi)$ satisfies the $d$--power condition, through the open $W$, then we have $d \mid N_j$ for all $j\in T$ satisfying $E_j \cap \overline{W} \neq \emptyset$.
\end{proposition}

\section{Description of $\nu$ via dual graph}

From now on we consider the plane curve case $n=2$, and we take $\pi$ as the {\em minimal} embedded resolution of $f$. In fact, we can as well study the germ of $f$ at the origin and allow $f$ to be an analytic function rather than just a polynomial.  At any rate, we slightly redefine $T$ as $T_e\cup T_s$, where $T_e$ runs over the {\em exceptional components} of $\pi$ and $T_s$ runs over the {\em analytically irreducible components of the strict transform} of $f$ by $\pi$.
In the (dual) resolution graph $\Gamma$ of $\pi$ one associates to each exceptional curve
$E_i$ a vertex, which we denote here for simplicity also by $E_i$,  and an arrowhead  to each (analytically)
irreducible component $E_i$ of the strict transform of $f$.
Each intersection between components $E_i$ is indicated by an edge
connecting the corresponding vertices or arrowhead.
We denote here by $\Gamma^e$ the restriction of $\Gamma$ to the exceptional locus, i.e., without the arrows.

For $i\in T_e$ we denote by $\delta_i$ the {\em valency} of $E_i$ in $\Gamma^e$, that is, the number of intersections of $E_i$ with other exceptional components, and hence also the number of edges in $\Gamma^e$ connected to $E_i$.


We use the language of Eisenbud--Neumann diagrams \cite{EN}, associated to the
(full) dual resolution graphs $\Gamma$ and $\Gamma^e$, where edges are decorated as follows.
For $i\in T_e$, an edge decoration $a$ 
next to $E_i$ along an edge $e$
adjacent to an exceptional $E_i$, indicates that $a$ is
the absolute value of the determinant of the intersection matrix of all exceptional components appearing in the subgraph of $\Gamma \setminus \{E_i\}$ in the direction of $e$.
These decorations satisfy the following properties.
\begin{itemize}
\item All edge decorations are positive integers.
\item The edge decorations along all edges next to a fixed $E_i$ are pairwise coprime and at most two of them are greater than $1$.
\item Fix an edge $e$ in $\Gamma$ between vertices $E_i$ and $E_j$ (thus corresponding to exceptional components). Let $a$ and $b$ be the decorations along $e$ next to
$E_i$ and $E_j$, respectively. Let also $a_k$ and $b_\ell$
denote the edge decorations along the other edges, connected to $E_i$ and
$E_j$, respectively. Then we have the {\em edge determinant rule}
$$ab - \prod_{k} a_k \prod_{\ell}b_\ell =1 $$
(where a product over the empty set is $1$).
\end{itemize}

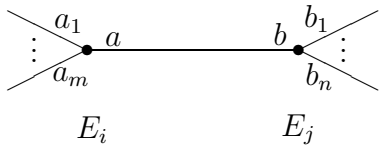
\begin{picture}(500,60)(50,-20)
\put(240,20){\circle*{4}}
\put(240,20){\line(1,0){80}}

\put(320,20){\circle*{4}}
\put(240,20){\line(-2, -1){30}}
 \put(240,20){\line(-2, 1){30}}
\put(220,23){\makebox(0,0){$\vdots$}}
\put(320,20){\line(2,-1){30}} \put(320,20){\line(2, 1){30}}
\put(337,23){\makebox(0,0){$\vdots$}}
\put(242,-10){\makebox(0,0){$E_i$}}
\put(320,-10){\makebox(0,0){$E_j$}}
\put(250,24){\makebox(0,0){$a$}} \put(313,26){\makebox(0,0){$b$}}
\put(327,32) {\makebox(0,0){$b_1$}}
\put(328,10){\makebox(0,0){$b_n$}}
 \put(233,30) {\makebox(0,0){$a_1$}}
 \put(234,10){\makebox(0,0){$a_m$}}
\end{picture}

\noindent
In fact, these properties are also valid for the dual graph of a non--minimal embedded resolution.  Using that $\pi$ is minimal we also have
\begin{itemize}
\item An edge decoration along an edge that is the start of a chain of exceptional components, ending in a vertex of valency $1$ of $\Gamma$, is greater than $1$.
\end{itemize}

\begin{example}\label{example}
Take $f = (y^2-x^3)^2-x^5y$. The decorated dual graph $\Gamma$ of its minimal embedded resolution is as follows, where we also indicate the numerical data $(N_i,\nu_i)$ of the $E_i$. Note that $c_0= 5/12$.

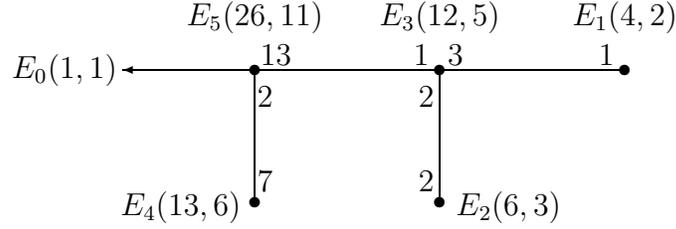
\begin{picture}(400,100)(-30,-35)

\put(150,25){\vector(-1,0){50}}
\put(150,25){\circle*{4}}
\put(220,25){\circle*{4}} \put(290,25){\circle*{4}}
\put(150,-25){\circle*{4}}
\put(220,-25){\circle*{4}}
\put(150,25){\line(1,0){140}}
\put(150,25){\line(0,-1){50}}
\put(220,25){\line(0,-1){50}}

\put(213,31){\makebox(0,0){$1$}}\put(158,31){\makebox(0,0){$13$}}
\put(154,15){\makebox(0,0){$2$}}
\put(154,-17){\makebox(0,0){$7$}}
\put(283,31){\makebox(0,0){$1$}}\put(226,31){\makebox(0,0){$3$}}
\put(215,-17){\makebox(0,0){$2$}}
\put(215,15){\makebox(0,0){$2$}}

\put(150,45){\makebox(0,0){$E_5(26,11)$}}
\put(220,45){\makebox(0,0){$E_3(12,5)$}}
\put(290,45){\makebox(0,0){$E_1(4,2)$}}
\put(122,-27){\makebox(0,0){$E_4(13,6)$}}
\put(246,-27){\makebox(0,0){$E_2(6,3)$}}
\put(78,25){\makebox(0,0){$E_0(1,1)$}}
\end{picture}
\end{example}

\medskip
We have the following well known \lq diagram calculus\rq, computing the numerical data $(N_i,\nu_i)$ of an exceptional curve $E_i$ in terms of the edge decorations of the graph $\Gamma$.
See for instance \cite{EN} and \cite{NV}. (It provides another way to compute the numerical data in Example \ref{example}.)

\begin{proposition}\label{N-theorem}
Fix an exceptional curve $E_i$. For any another component $E_j$, let $\ell_{ij}$ be the product of the edge decorations that are adjacent to, but not on, the path in $\Gamma$ from $E_i$ to $E_j$. Then
\begin{equation}\label{N-formula}
N_i = \sum_{j\in T_s} \ell_{ij} N_j ,
\end{equation}
\begin{equation}\label{old-nu-formula}
\nu_i = \sum_{j\in T_e} \ell_{ij} (2-\delta_j).
\end{equation}
\end{proposition}

\smallskip
We now show a useful upper bound for $\nu_i$, depending only on the edge decorations along $E_i$, that is often even an equality.

\begin{theorem}\label{nu-theorem}
(1) Let $E$ be a vertex of valency at least $2$ in $\Gamma^e$. Say $a$ and $b$ are edge decorations at $E$ such that all other edge decorations at $E$ are $1$. (Possibly also $a$ or $b$ are $1$.) Then we have $\nu \leq a+b$. More precisely, we have the following.

(i) If starting from $E$, say in the $a$--decorated edge direction, there exists in some part of $\Gamma^e$ a vertex of valency at least $3$ as in the figure below, where both $c>1$ and $d>1$,  then $\nu \leq a-b$.

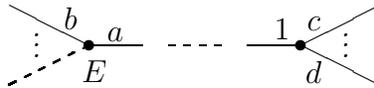
\begin{picture}(500,60)(50,-10)
\put(240,20){\circle*{4}}
\dashline[3]{3}(270,20)(290,20)
\put(320,20){\circle*{4}}
\put(240,20){\line(1,0){20}}\put(300,20){\line(1,0){20}}
\dashline[3]{3}(240,20)(210,5)
 \put(240,20){\line(-2, 1){30}}
\put(220,23){\makebox(0,0){$\vdots$}}
\put(320,20){\line(2,-1){30}} \put(320,20){\line(2, 1){30}}
\put(337,23){\makebox(0,0){$\vdots$}}
\put(242,10){\makebox(0,0){$E$}}
\put(250,24){\makebox(0,0){$a$}} \put(313,26){\makebox(0,0){$1$}}
\put(325,28) {\makebox(0,0){$c$}}
\put(325,10){\makebox(0,0){$d$}}
 \put(233,30) {\makebox(0,0){$b$}}
\end{picture}

(ii) If there is no such vertex (in any direction starting from $E$), then $\nu=a+b$.

\smallskip
\noindent
(2) Let $E$ be a vertex of valency $1$ in $\Gamma^e$, with edge decoration $a$. Then $\nu \leq a+1$, and more precisely, with the analogous case distinction, either $\nu \leq a-1$ or $\nu = a+1$.
\end{theorem}

\noindent
{\em Note.} There can be at most one direction as in (i) starting from $E$, which is well known and will also be clear from the proof.

\begin{proof}
Consider any vertex $E_j$ of $\Gamma^e$
 with some adjacent edge decoration equal to $1$, such that the subgraph $\Gamma_j$ in the direction of this edge does {\em not} contain $E$. (Possibly $E_j=E$.)

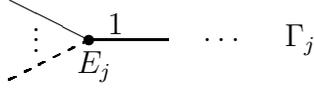
\begin{picture}(450,60)(50,-10)
\put(240,20){\circle*{4}}
\dashline[3]{3}(240,20)(210,5)
\put(320,20){\makebox(0,0){$\Gamma_j$}}
\put(240,20){\line(1,0){30}}
\put(290,20){\makebox(0,0){$\hdots$}}
 \put(240,20){\line(-2, 1){30}}
\put(220,23){\makebox(0,0){$\vdots$}}
\put(242,10){\makebox(0,0){$E_j$}}
\put(250,26){\makebox(0,0){$1$}}
\end{picture}

\noindent
We claim that, in order to compute $\nu$, we can contract/forget the subgraph $\Gamma_j$.
Indeed, since the absolute value of the determinant of the intersection matrix of $\Gamma_j$ is $1$, all the exceptional curves in  $\Gamma_j$ can be blown down. We can consider this \lq blown down situation\rq\ as an intermediate step in constructing $\pi$, and $\nu$ can be computed on the graph of that intermediate step. (Alternatively, one can prove the claim using an elementary computation with formula (\ref{old-nu-formula}).)

Now we contract/delete all such subgraphs. The resulting graph $\Gamma_0$, corresponding to some intermediate step in constructing $\pi$, must satisfy one of the two following properties.

\smallskip
(i) There is still a vertex of valency at least $3$ in $\Gamma_0$, say in the $a$--decorated edge direction. Then $\Gamma_0$ is necessarily of the form below, where all $a_i, b_i >1$. When $b>1$, the part of $\Gamma_0$ in the $b$--decorated edge direction is a chain, and when $b=1$, the vertex $E$ has valency $1$ in $\Gamma_0$.

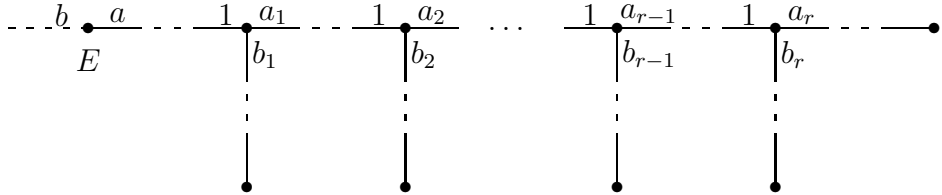
\begin{picture}(405,85)(-5,-30)
 \put(50,40){\circle*{4}}
\put(110,40){\circle*{4}} \put(170,40){\circle*{4}}
\put(250,40){\circle*{4}} \put(310,40){\circle*{4}}
\put(370,40){\circle*{4}}
 \put(110,-20){\circle*{4}}
\put(170,-20){\circle*{4}} \put(250,-20){\circle*{4}}
\put(310,-20){\circle*{4}}


\dashline[3]{3}(20,40)(50,40)
\put(50,40){\line(1,0){20}}
\dashline[3]{3}(75,40)(85,40)
\put(90,40){\line(1,0){40}}
\dashline[3]{3}(135,40)(145,40)
\put(150,40){\line(1,0){40}}
\put(230,40){\line(1,0){40}}
\dashline[3]{3}(275,40)(285,40)
\put(290,40){\line(1,0){40}}
\dashline[3]{3}(335,40)(345,40)
\put(350,40){\line(1,0){20}}

\put(170,40){\line(0,-1){20}}
\dashline[3]{3}(170,15)(170,5)
\put(170,0){\line(0,-1){20}}

\put(110,40){\line(0,-1){20}}
\dashline[3]{3}(110,15)(110,5)
\put(110,0){\line(0,-1){20}}

\put(250,40){\line(0,-1){20}}
\dashline[3]{3}(250,15)(250,5)
\put(250,0){\line(0,-1){20}}

\put(310,40){\line(0,-1){20}}
\dashline[3]{3}(310,15)(310,5)
\put(310,0){\line(0,-1){20}}

\put(40,45){\makebox(0,0){$b$}}
\put(102,45){\makebox(0,0){$1$}}
\put(160,45){\makebox(0,0){$1$}}
\put(240,45){\makebox(0,0){$1$}}
\put(300,45){\makebox(0,0){$1$}}
\put(120,45){\makebox(0,0){$a_1$}}
\put(180,45){\makebox(0,0){$a_2$}}
\put(262,45){\makebox(0,0){$a_{r-1}$}}
\put(320,45){\makebox(0,0){$a_r$}}

\put(58,45){\makebox(0,0)[l]{$a$}}
\put(112,30){\makebox(0,0)[l]{$b_1$}}
\put(172,30){\makebox(0,0)[l]{$b_2$}}
\put(252,30){\makebox(0,0)[l]{$b_{r-1}$}}
\put(312,30){\makebox(0,0)[l]{$b_r$}}

\put(208,40){\makebox(0,0){$\ldots$}}

\put(50,28){\makebox(0,0){$E$}}
\end{picture}

\noindent
Note that in our resolution graphs there can be at most one vertex of valency at least $3$ with two attached chains and both edge decorations larger than $1$
(as on the most right part of the figure above). Indeed, by a contraction argument as before, we can consider the subgraph consisting of only that vertex and the two attached chains as corresponding to some intermediate step of $\pi$, and then this subgraph must contain the first created exceptional curve as a vertex.

Using formula (\ref{old-nu-formula}) we have
$$\aligned
\nu = \, &a+b[ a_1-a_1b_1+(a_2-a_2b_2)b_1+\dots+(a_{r-1}-a_{r-1}b_{r-1})b_1b_2\dots b_{r-2} \\
&+ (a_r+b_r-a_rb_r)b_1b_2\dots b_{r-1}].
\endaligned$$
Since all $a_i-a_ib_i$ and also $a_r+b_r-a_rb_r$ are negative, we have that $\nu \leq a-b$.

\smallskip
(ii) There is no vertex of valency $3$ in $\Gamma_0$. Then we have $\nu=a+b$. (If $a$ or $b$ are equal to $1$, then $E$ has valency at most $1$ in $\Gamma_0$.)

\smallskip\noindent
The proof of (2) is completely analogous.

\end{proof}

\smallskip
\begin{example}[continuing Example \ref{example}]
All different cases of Theorem \ref{nu-theorem} occur in the example. In particular, $E_5$ and $E_4$ satisfy the inequality \lq $\nu\leq a-b$\rq, and both inequalities are sharp here.
\end{example}

\section{Proof of the main theorem}

In dimension 2 we only have the cases $|I|=1$ and $|I|=2$. By Remark \ref{trivial} and Proposition \ref{d-divisible} we may and will assume that
\begin{itemize}
\item  when $I=\{i\}$, the component $E_i$ is exeptional and does not intersect the strict transform of $f$,
\item when $I=\{i,j\}$ (with $i\neq j$), the components $E_i$ and $E_j$ are exceptional.
\end{itemize}

\noindent
We will in fact show a slightly stronger statement than (\ref{formulation}).

\begin{theorem}\label{main theorem}
Let $d \in \Z_{\geq 2}$.

(1) Let $E_i$ be an exceptional component such that $d \mid N_i$ and  $d \mid N_\ell$ for all components $E_\ell$ intersecting $E_i$.
Then either
$$\frac{\nu_i}{N_i} \leq \frac 1d \qquad\text{or}\qquad \frac {\nu_\ell}{N_\ell} \leq \frac{ \nu_i + 1/d}{ N_i + 1}$$
for some intersecting component $E_\ell$.

(2) Let $E_i$ and $E_j$ be intersecting exceptional components such that $d \mid N_i$ and  $d \mid N_j$.
Then (up to a switch of the indices)
$$\frac{\nu_i}{N_i}\leq \frac 1d \qquad\text{or}\qquad \frac {\nu_j}{N_j} \leq \frac{ \nu_i + 1/d}{ N_i + 1} .$$
\end{theorem}

\bigskip
More precisely, we will argue by case distinction, depending on the position of $E_i$ and $E_j$ in the graph $\Gamma$. These different results could be of interest for future reference; for that reason we formulate them in separate independent statements.

\begin{lemma}\label{lemma1}
Let $E_i$ and $E_j$ be adjacent vertices on the graph $\Gamma$. Let $d \in \Z_{\geq 2}$ such that $d \mid N_i$ and  $d \mid N_j$.
Suppose that there exist arrows in $\Gamma$ on both sides of the edge between $E_i$ and $E_j$. Then $\nu_i/N_i \leq 1/d $ and $\nu_j/N_j \leq 1/d $.
\end{lemma}

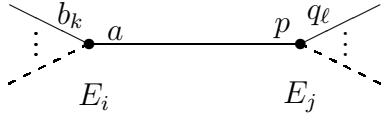
\begin{picture}(500,60)(50,-10)
\put(220,23){\makebox(0,0){$\vdots$}}
\put(240,20){\circle*{4}}
\dashline[3]{3}(320,20)(350,5)
\put(320,20){\circle*{4}}
\put(240,20){\line(1,0){80}}
\dashline[3]{3}(240,20)(210,5)
 \put(240,20){\line(-2, 1){30}}
\put(220,23){\makebox(0,0){$\vdots$}}
 \put(320,20){\line(2, 1){30}}
\put(337,23){\makebox(0,0){$\vdots$}}
\put(242,0){\makebox(0,0){$E_i$}}
\put(320,0){\makebox(0,0){$E_j$}}
\put(250,24){\makebox(0,0){$a$}} \put(313,26){\makebox(0,0){$p$}}
\put(327,31) {\makebox(0,0){$q_\ell$}}
 \put(233,31) {\makebox(0,0){$b_k$}}
\end{picture}

\begin{proof} Let $a$ and $p$ be the edge decorations at $E_i$ and $E_j$ on the edge connecting them, and $b_k$ and $q_\ell$ the other edge decorations at $E_i$ and $E_j$, respectively.
By Proposition \ref{N-theorem} we have that
$$N_i = La+R\prod_k b_k \qquad\text{and}\qquad N_j=L\prod_\ell q_\ell+Rp ,$$
where $L$ and $R$ describe the total contribution in formula (\ref{N-formula}) of arrows \lq on the left of $E_i$\rq\ and \lq on the right of $E_j$\rq, respectively.
Since $ap-\prod_k b_k \prod_\ell q_\ell=1$ (edge determinant rule), we derive that
$$N_i p - N_j \prod_k b_k  = L \qquad\text{and}\qquad N_j a - N_i \prod_\ell q_\ell =R ,$$
and hence that $d \mid L$ and $d \mid R$.

If two of the $b_k$ are greater than $1$, say $b_1$ and $b_2$, then  $a=1$ and $N_i=L+Rb_1b_2$. By Theorem \ref{nu-theorem} we have that $\nu_i \leq b_1+b_2$ and
consequently
$$ \frac{\nu_i}{N_i} \leq \frac{b_1+b_2}{d(\frac Ld  + \frac Rd b_1b_2)} < \frac 1d .$$
If on the other hand at most one of the $b_k$ is greater than $1$, say (at most) $b_1$, then $N_i=La+Rb_1$. Now we have by Theorem \ref{nu-theorem}  that $\nu_i \leq a+b_1$ and
 consequently
$$ \frac{\nu_i}{N_i} \leq \frac{a+b_1}{d(\frac Ld a + \frac Rd b_1)} \leq \frac 1d .$$
By symmetry the same results holds for $\nu_j/N_j$.
\end{proof}

\smallskip
\begin{example} Take $f=x^2(y^2-x^4)$. Its minimal embedded resolution provides an easy illustration of Lemma \ref{lemma1} with $d=2$, where moreover both inequalities are sharp.

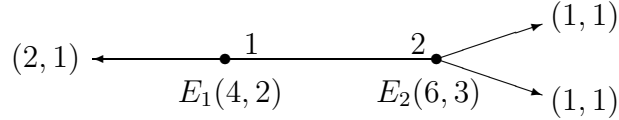
\begin{picture}(400,60)(50,-10)
\put(240,20){\circle*{4}}
\put(320,20){\circle*{4}}
\put(240,20){\line(1,0){80}}
\put(240,20){\vector(-1,0){50}}
 \put(320,20){\vector(3, 1){40}}
 \put(320,20){\vector(3, -1){40}}

\put(376,36){\makebox(0,0){$(1,1)$}}
\put(376,4){\makebox(0,0){$(1,1)$}}
\put(172,20){\makebox(0,0){$(2,1)$}}
\put(242,7){\makebox(0,0){$E_1(4,2)$}}
\put(317,7){\makebox(0,0){$E_2(6,3)$}}
\put(250,26){\makebox(0,0){$1$}} \put(313,26){\makebox(0,0){$2$}}

\end{picture}

\end{example}

\smallskip
\begin{lemma}\label{lemma2}
Let $E_i$ and $E_j$ be adjacent vertices on the graph $\Gamma$. Let $d \in \Z_{\geq 2}$ such that $d \mid N_i$ and  $d \mid N_j$.
Suppose that, besides the edge in the direction of $E_i$, the vertex $E_j$ is adjacent precisely to a subgraph of $\Gamma$ of the following form (where the valency of $E_j$ in $\Gamma$ can be $2$ or $3$, and there is at least one vertical chain). Then $\nu_i/N_i < 1/d $.
\end{lemma}

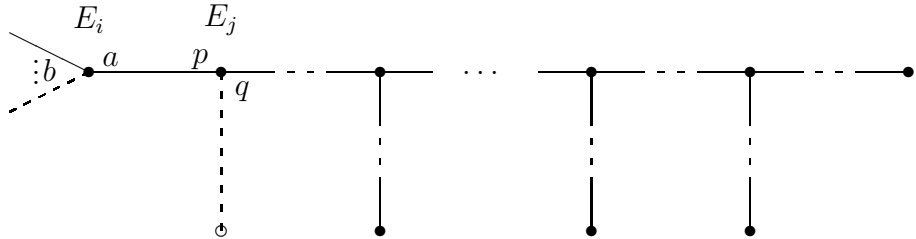
\begin{picture}(230,110)(-15,-30)
\put(40,43){\makebox(0,0){$\vdots$}}
 \put(60,40){\circle*{4}}
\put(110,40){\circle*{4}} \put(170,40){\circle*{4}}
\put(250,40){\circle*{4}} \put(310,40){\circle*{4}}
\put(370,40){\circle*{4}}
 \put(110,-20){\circle{4}}
\put(170,-20){\circle*{4}} \put(250,-20){\circle*{4}}
\put(310,-20){\circle*{4}}

\dashline[3]{3}(60,40)(30,25)
 \put(60,40){\line(-2, 1){30}}

\put(60,40){\line(1,0){70}}
\dashline[3]{3}(135,40)(145,40)
\put(150,40){\line(1,0){40}}
\put(230,40){\line(1,0){40}}
\dashline[3]{3}(275,40)(285,40)
\put(290,40){\line(1,0){40}}
\dashline[3]{3}(335,40)(345,40)
\put(350,40){\line(1,0){20}}

\put(170,40){\line(0,-1){20}}
\dashline[3]{3}(170,15)(170,5)
\put(170,0){\line(0,-1){20}}

\dashline[3]{3}(110,40)(110,-20)

\put(250,40){\line(0,-1){20}}
\dashline[3]{3}(250,15)(250,5)
\put(250,0){\line(0,-1){20}}

\put(310,40){\line(0,-1){20}}
\dashline[3]{3}(310,15)(310,5)
\put(310,0){\line(0,-1){20}}

\put(45,41){\makebox(0,0){$b$}}
\put(102,45){\makebox(0,0){$p$}}

\put(118,33){\makebox(0,0){$q$}}
\put(65,45){\makebox(0,0)[l]{$a$}}

\put(208,40){\makebox(0,0){$\ldots$}}

\put(60,60){\makebox(0,0){$E_i$}}
\put(110,60){\makebox(0,0){$E_j$}}
\end{picture}

\begin{proof} Let $q$ be the decoration or the product of the two decorations at $E_j$, not on the edge between $E_i$ and $E_j$.
By Proposition \ref{N-theorem} we have that
$N_i = La$ and $N_j=Lq$,
where $L$ is the total contribution in formula (\ref{N-formula}) of all arrows in $\Gamma$.
Since $ap-bq=1$, we derive that
$pN_i -bN_j = L$ and hence that $d \mid L$.  By Theorem \ref{nu-theorem} we have that $\nu_i \leq a-b$ and consequently
$$\frac{\nu_i}{N_i}  \leq \frac{a-b}{La} < \frac a{La} = \frac 1L   \leq \frac 1d .$$
\end{proof}

\smallskip
\begin{example}[continuing Example \ref{example}]
The vertices $E_5$ and $E_3$ form an illustration of Lemma \ref{lemma2} with $d=2$.
\end{example}

\smallskip
\begin{lemma}\label{lemma3}
Let $E_1$ be an end vertex of $\Gamma$, such that $E_1, E_2, \dots, E_r$ form a chain in $\Gamma$ (with $r\geq 2$). Let $d \in \Z_{\geq 2}$ such that $d \mid N_r$ and  $d \mid N_{r-1}$.
Then
\begin{equation}\label{inequality3}
\frac{\nu_r}{N_r} \leq \frac{ \nu_{r-1} + 1/d}{ N_{r-1} + 1} .
\end{equation}
\end{lemma}

\begin{picture}(230,85)(-45,0)
\put(40,43){\makebox(0,0){$\vdots$}}
 \put(60,40){\circle*{4}}
\put(115,40){\circle*{4}}
 \put(215,40){\circle*{4}}
\put(270,40){\circle*{4}}

\dashline[3]{3}(60,40)(30,25)
 \put(60,40){\line(-2, 1){30}}

\put(60,40){\line(1,0){80}}
\put(190,40){\line(1,0){80}}

\put(54,51){\makebox(0,0){$b_r$}}
\put(101,46){\makebox(0,0){$b_{r-1}$}}
\put(128,45){\makebox(0,0){$a_{r-1}$}}
\put(65,45){\makebox(0,0)[l]{$a_r$}}

\put(205,46){\makebox(0,0){$b_2$}}
\put(260,46){\makebox(0,0){$b_{1}$}}
\put(220,45){\makebox(0,0)[l]{$a_2$}}

\put(165,40){\makebox(0,0){$\ldots$}}

\put(60,20){\makebox(0,0){$E_r$}}
\put(115,20){\makebox(0,0){$E_{r-1}$}}

\put(215,20){\makebox(0,0){$E_2$}}
\put(270,20){\makebox(0,0){$E_{1}$}}

\end{picture}

\noindent
{\em Note.} The condition $d \mid N_r$ and  $d \mid N_{r-1}$ is equivalent to for instance $d \mid N_1$ and to $d \mid N_i$ for all $i=1\dots,r$.

\begin{proof}
It is well known (see e.g. \cite{V}) that $N_r= a_rN_1$ and $N_{r-1} = a_{r-1}N_1$ (it also follows from Proposition \ref{N-theorem}).  By Theorem \ref{nu-theorem} we have
$$ \nu_r = b_r +xa_r \qquad\text{and}\qquad \nu_{r-1} = b_{r-1} +xa_{r-1} ,$$
where $x=1$ or $x \in \Z_{<0}$, in particular $x \leq 1$.
Substituting these equalities in (\ref{inequality3}) yields, after a straightforward calculation,
$$
b_r \leq N_1(a_rb_{r-1} - b_ra_{r-1}) + (\frac{N_1}d -x) a_r = N_1 + (\frac{N_1}d -x) a_r .
$$
Since $x \leq 1$ and $d \mid N_1$, this will be implied by
\begin{equation}\label{inequality}
b_r \leq N_1 .
\end{equation}
Say the chain above ends in a vertex $E_n$ of valency at least $3$ in $\Gamma$ (where $n\geq r$). It is well known and easily verified that
$ a_{i+1} \geq a_i +1$ for $i=1,\dots,n-1$, where $a_1=1$ (it follows from the fact that, in the minimal embedded resolution, all self--intersections $E_1^2,\dots, E_{n-1}^2$ are at most $-2$).
Then, by an elementary calculation, this implies $b_{i} \leq b_{i+1} $ for $i=1,\dots,n-1$ and in particular $b_{r} \leq b_n$.

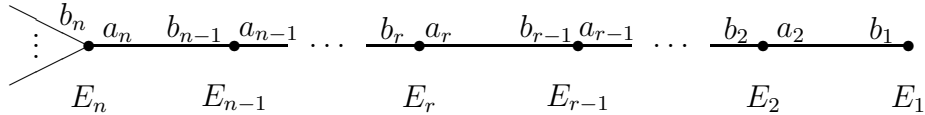
\begin{picture}(230,80)(-15,0)
\put(40,43){\makebox(0,0){$\vdots$}}
 \put(60,40){\circle*{4}}
\put(115,40){\circle*{4}} \put(185,40){\circle*{4}}
\put(245,40){\circle*{4}} \put(315,40){\circle*{4}}
\put(370,40){\circle*{4}}

\put(60,40){\line(-2,-1){30}}
\put(60,40){\line(-2,1){30}}

\put(60,40){\line(1,0){75}}
\put(165,40){\line(1,0){100}}
\put(295,40){\line(1,0){75}}

\put(54,51){\makebox(0,0){$b_n$}}
\put(100,46){\makebox(0,0){$b_{n-1}$}}

\put(128,45){\makebox(0,0){$a_{n-1}$}}
\put(65,45){\makebox(0,0)[l]{$a_n$}}

\put(175,46){\makebox(0,0){$b_r$}}
\put(233,46){\makebox(0,0){$b_{r-1}$}}
\put(256,45){\makebox(0,0){$a_{r-1}$}}
\put(187,45){\makebox(0,0)[l]{$a_r$}}

\put(305,46){\makebox(0,0){$b_2$}}
\put(360,46){\makebox(0,0){$b_{1}$}}
\put(320,45){\makebox(0,0)[l]{$a_2$}}

\put(150,40){\makebox(0,0){$\ldots$}}
\put(280,40){\makebox(0,0){$\ldots$}}

\put(60,20){\makebox(0,0){$E_n$}}
\put(115,20){\makebox(0,0){$E_{n-1}$}}
\put(185,20){\makebox(0,0){$E_r$}}
\put(245,20){\makebox(0,0){$E_{r-1}$}}
\put(315,20){\makebox(0,0){$E_2$}}
\put(370,20){\makebox(0,0){$E_{1}$}}

\end{picture}

We claim that $b_n \leq N_1$, which then implies (\ref{inequality}). When $b_n=1$, this is trivial. Otherwise the decorations along the other edges adjacent to $E_n$ are $1$, implying that, in the direction of such an other edge away from $E_n$, there is at least one arrow. And then Proposition \ref{N-theorem} yields that $N_1 \geq b_n$.
\end{proof}

\smallskip
\begin{example}[continuing Example \ref{example}]
The end vertices $E_1$ with $d=4$, $E_2$ with $d=6$, and $E_4$ with $d=13$, form three different illustrations of Lemma \ref{lemma3}. Each time the length $r$ of the chain is just $2$. Here the inequalities are not sharp. (The minimal embedded resolution of $f=y^2-x^3$ provides two sharp examples.)
\end{example}

\smallskip
\begin{lemma}\label{lemma4}
Let $E_i$ be a vertex of the graph $\Gamma$ of valency at least $3$, where two chains are attached to $E_i$ with end vertices $E_1$ and $E_2$, respectively.
Let $d \in \Z_{\geq 2}$ such that $d \mid N_1$ and  $d \mid N_2$, and hence $d \mid N_i$.
 Then $\nu_i/N_i < 1/d $.
\end{lemma}

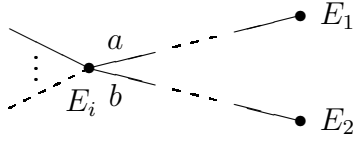
\begin{picture}(500,65)(50,-20)
\put(220,23){\makebox(0,0){$\vdots$}}
\put(240,20){\circle*{4}}
\dashline[3]{3}(240,20)(210,5)
 \put(240,20){\line(-2, 1){30}}
 \put(220,23){\makebox(0,0){$\vdots$}}

\put(240,20){\line(4,1){25}}
\dashline[3]{3}(272,28)(288,32)
\put(320,40){\line(-4,-1){25}}
 \put(240,20){\line(4, -1){25}}
\dashline[3]{3}(272,12)(288,8)
\put(320,0){\line(-4,1){25}}

\put(320,40){\circle*{4}}
\put(320,0){\circle*{4}}

\put(237,7){\makebox(0,0){$E_i$}}
\put(334,40){\makebox(0,0){$E_1$}}
\put(334,0){\makebox(0,0){$E_2$}}

\put(250,30){\makebox(0,0){$a$}}

 \put(250,10) {\makebox(0,0){$b$}}
\end{picture}

\noindent
{\em Note.} The conditions $d \mid N_1$ and $d \mid N_2$ imply that  $d \mid N_\ell$ for all $E_\ell$ in both chains.

\begin{proof} By Proposition \ref{N-theorem} we have that
$$N_i = Lab, \quad N_1=Lb, \quad N_2=La,$$
where $L$ is the total contribution in formula (\ref{N-formula}) of the arrows \lq on the left of $E_i$\rq.
Since $a$ and $b$ are coprime, we derive that
$d \mid L$.  By Theorem \ref{nu-theorem}, we have that $\nu_i = a +b$. Hence
$$ \frac{\nu_i}{N_i} = \frac{a+b}{Lab} <  \frac 1L   \leq \frac 1d ,$$
where the inequality follows from $a,b > 1$.
\end{proof}

\smallskip
\begin{example}[continuing Example \ref{example}]
The vertex $E_3$ has two such attached chains and illustrates Lemma \ref{lemma4} with $d=2$.
\end{example}

\medskip
We now show Theorem \ref{main theorem}.

\begin{proof}
(1) If $E_i$ satisfies the condition of Lemma \ref{lemma1}, then $\nu_i/N_i \leq 1/d$. Otherwise, there is only one \lq arrow direction\rq\ starting from $E_i$.

When $E_i$ has valency $1$ (in $\Gamma$), we can consider it as the vertex $E_1$ in Lemma \ref{lemma3}, and then its adjacent vertex $E_2$ satisfies
$$\frac{\nu_2}{N_2} \leq \frac{ \nu_{1} + 1/d}{ N_{1} + 1} .$$
When $E_i$ has valency at least $2$, we can either consider it as the vertex $E_{r-1}$ of Lemma \ref{lemma3}, or 
 it must be in the situation of Lemma \ref{lemma2} or Lemma \ref{lemma4}. In all these cases the conclusion holds.

(2) If $E_i$ and $E_j$ satisfy the condition of Lemma \ref{lemma1}, both $\nu_i/N_i \leq 1/d$ and $\nu_j/N_j \leq 1/d$. Otherwise they satisfy either the condition of Lemma \ref{lemma2} or they can be identified with the vertices $E_r$ and $E_{r-1}$ in the chain of Lemma \ref{lemma3}, implying one of the desired inequalities.
\end{proof}

\end{document}